# A branch and bound algorithm for the robust parallel machine scheduling problem with sequence-dependent setup time

Lei Liu, Canrong Zhang

Research Center for Modern Logistics, Graduate School at Shenzhen, Tsinghua University

**Abstract:** This paper focuses on the identical parallel machine scheduling problem with sequence-dependent setup time, with special attention paid to the uncertainty of processing time. In this paper, a mathematical model of the parallel machine scheduling problem with stochastic processing time is constructed. Then a branch and bound (B&B) algorithm is proposed, and three methods for generating upper bounds and one heuristic for generating lower bound are developed to evaluate the performance of the B&B. In addition, a new branch scheme and two dominance rules are also devised to further improve the efficiency of the algorithm. In the numerical experiment analysis, the performance of the new branch scheme, the upper and lower bounds, and the dominance rule are tested, and the robustness of the stochastic parallel machine scheduling problem is verified by compared it with the deterministic version .

## 1. Introduction

In the semi-conductor factory, parallel machine scheduling with sequence dependent setup time is frequently encountered. Great research efforts are paid to this filed in the research community [1-3]. Most studies on parallel machine scheduling problems, however, are limited to deterministic scheduling despite the environment full of uncertainty. Just due to the uncertainties such as the processing-time uncertainty, the deterministic scheduling using nominal value as input likely loses its appealing because the uncertainty has been compromised. Therefore, a new method is needed to tackle uncertainty explicitly.

In this paper, we examine the identical parallel machine scheduling with sequence dependent setup time and stochastic processing time, with a goal to find a robust schedule that maximizes the service level for customers, which is represented as the product of the probabilities that the completion time of each customer does not exceed their respective pre-specified makespans.

The contributions of this paper in the following three points: (1) a mathematic model is constructed for the parallel machine scheduling with sequence dependent setup time and stochastic processing time; (2) an exact B&B algorithm is proposed to solve the problem; and (3) new branch scheme, upper and lower bounds, as well as dominance rules are derived to accelerate convergence of the B&B.

The remainder of this paper is organized as follows. Section 2 gives a literature review of scheduling problems under stochastic situations; Section 3,provides a description of the problem and a non-linear binary formulation; Section 4 describes the branch and bound algorithm as well as the branch scheme, low bound and upper bounds construction, and dominance rules. Section 5 presents the computational experiments and some findings which are drawn from the numerical results; Section 6 provides conclusions and suggestions for future work.

## 2. Literature Review

Multiple articles have been published on the stochastic machine scheduling problem with objectives such as minimization of makespan, minimization of total flow time and the due date related objectives, the methods used in these articles include fuzzy method, robust method and so on[4-5].

The "$\beta$-robust" scheduling model was first introduced in the single machine environment, Daniels and Kouvelis[6] assumed uncertain processing time for jobs in a single machine environment and developed a robust scheduling solution to single machine problem using the total flow time as the objective function. They developed a branch and bound algorithm and two heuristic algorithms to achieve robust schedule. Gerodimos[7] defined three robustness measures for the scheduling problem and developed a mathematical model for the problem of single machine scheduling and a two-machine flow shop. Wu, Brown and Beck[8] proposed a robust scheduling model with constraint programming for the single machine problem, they developed three constraint programming models and solved this problem up to 20 jobs. In recent years, Chung-Cheng Lu[9] and Paolo Detti[10] also used the robust methods to solve the single machine scheduling problem with uncertainty.

In the robust parallel machine scheduling area, Ranjbar, Davari and Leus[11] used the $\beta$-robust scheduling approach to solve the scheduling problem of identical parallel machines with uncertain processing times. That study found a robust scheduling that maximize the customer service level. Two branch and bound algorithms were developed for finding an optimal solution. S. Alimoradi and M. Hematian[12] studied the problem of maximizing the probability that the total flow time does not exceed a predetermined limit in identical parallel machines while the processing time of each job is stochastic. They proposed several theorems and the theorems considerably reduced the search space and resulted in a branch and bound method. They also proposed several dominance rules, an upper bound, and a lower bound are proposed to help find optimal solutions to the problems. Mir Saber Salehi and Javad Rezasian[13] used a meta heuristic (PSO-GA) approach to solve the robust unrelated parallel machine scheduling problem.

## 3. Problem description and formulation

Based on the model put forward by M.Ranjbar[9], our model takes into account the sequence-dependent setup time, making the model better fit the practical situation. Because of the sequence dependent setup time, the order of the jobs assigned to the same machine will significant impact on the objective value. This is the major difference between this paper and [11-12].

In the problem, there are n independent jobs represented by set $J=\{J_1, J_2, ..., J_n\}$, and a set of m identical machines $M=\{M_1, M_2, ..., M_m\}$. Each job has a stochastic processing time $p_j$ following a normal distribution function $f(p_j)$, with mean $E(p_j) = \mu_j$ and variance $Var(p_j) = \sigma_j^2$ $(j=1,...,n)$. Sequence dependent setup time $s_{ij}$ represents the time for setting up the machine for job j from the state processing job i. Normally, $s_{ij}$ is not equal to $s_{ji}$. During the processing, preemption is not allowed. We assume that all the jobs are released at the starting of the planning horizon, and a uniform due time is set for all the jobs. the customer service level $P$ is set as the product of the probabilities that the due day is not exceeded by the completion times of customers. The objective of our model is to maximize the customer service level.

We represent the total setup time on machine $M_i$ by $S_i$, the total time $T_i$ on $M_i$ has a normal distribution: $T_i \sim N(\mu_{M_i} + S_i, \sigma_{M_i}^2)$, and the customer service level on $M_i$ can be written as $\pi_i = P(T_i \leq \delta)$.

Then we use $AD\mu_{Mi} = \mu_{Mi} + S_i$ represent the expectation of time on machine $M_i$, and change the normal distribution to a standard normal distribution: $Z_i = (T_i - AD\mu_{Mi})/\sigma_{Mi}$, then the customer service level on $M_i$ becomes $\pi_i = P(Z_i \leq z_i) = \phi(z_i)$.

Then we introduce the following decision variables:

$$X_{ij} = \begin{cases} 1, & \text{if } J_j \text{ is assigned to } M_i \\ 0, & \text{otherwise} \end{cases} \quad i=1,\ldots,n; \, j=1,\ldots m$$

$$Y_{jk} = \begin{cases} 1, & \text{if } J_j \text{ is proceeded after } J_k \text{ on the same machine} \\ 0, & \text{otherwise} \end{cases} \quad j=1,\ldots,n; \, k=1,\ldots,n; \, j \neq k$$

We have $P = P(max\{T_1,\ldots,T_m\} \leq \delta)$, and $T_1, T_2, \ldots, T_m$ are independent variables, we can rewrite P as $P = P(T_1 \leq \delta)P(T_2 \leq \delta)\ldots P(T_m \leq \delta)$, and the objective function can be rewrite as $\prod_{i=1}^{m} \pi_i$, which maximize the customer service level as well as balance the machine workload. Then the robust parallel machine scheduling problem with sequence dependent setup time can now be formulated as the following non-linear binary model:

$$MaxP = \prod_{i=1}^{m} \phi(\frac{\delta - AD\mu_{Mi}}{\sigma_{Mi}}) \tag{1}$$

S.t

$$\sum_{i=1}^{m} X_{ij} = 1, \quad j = 1, 2, \ldots, n \tag{2}$$

$$AD\mu_{Mi} = \sum_{j=1}^{n} \mu_j X_{ij} + \sum_{j=1}^{n} S_{jk} Y_{jk}, \quad i = 1, 2, \ldots, m \tag{3}$$

$$\sigma_{Mi}^2 = \sum_{j=1}^{n} \sigma_j^2 X_{ij} \quad i = 1, 2, \ldots, m \tag{4}$$

$$X_{ij} \in \{0,1\} \quad i = 1, \ldots, m; j = 1, \ldots, n \tag{5}$$

$$Y_{jk} \in \{0,1\} \quad j = 1, \ldots, n; k = 1, \ldots, n \tag{6}$$

Constraint (2) assures that each job is assigned to exactly one machine. Constraints (3) and (4) compute the mean and the variance of total time on each machine, which are used to determine P. Constraint (5) and Constraint (6) states that all the decision variables are binary.

## 4. The branch and bound method

### *4.1 The branch scheme*

We develop a best-first branch and bound algorithm to find the optimal solution of the problem. Because of the setup time, we cannot use the branch scheme introduced in the literature without setup time like M.Ranjbar's work[9] .

In the B&B trees for most parallel machine scheduling problem with sequence dependent setup time, each node corresponds to a partial schedule, and a node in the *i*th level of the tree represents a partial schedule in which i jobs are scheduled at the front part of a complete schedule. In the *i*th level, m·(n-i) child nodes are generated from the parent node, where m and n are the number of machines and the number of jobs, respectively. We use a two machines and two jobs problem as an example in figure 1.

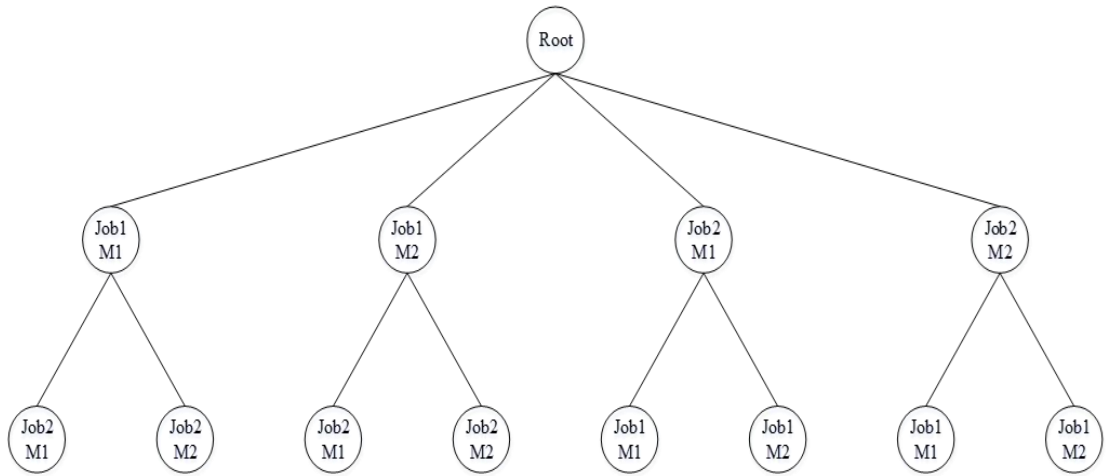

Figure 1

In this section, we suggest a new branch scheme for the robust parallel machine scheduling problem. Firstly, we use the branch scheme introduced by S.O. Shim and Y.D. Kim[12] to reduce the redundant nodes. In this scheme, a new job (associated with a certain node in the B&B tree) is assigned to a machine whose index is not less than that of the machine associated with the parent node. In other words, from each parent node, we generate child nodes associated with the parent node, and do not generate child nodes that do not satisfy this condition. In a two jobs and two machines example, we can generate a tree in figure 2.

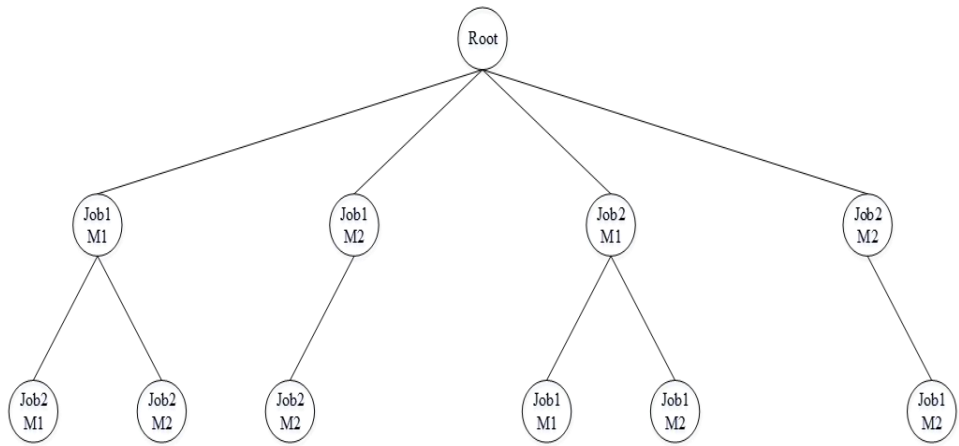

Figure 2

Farther, we can see that on the *i*th level, if the index of the assigned machine is bigger than *i*, it means that in this branch tree we will have at least one idle machine at last, because of the objective function, if one of the machines is idle, it must not be the optimal schedule. Based on this, we can simplify the example branch tree as figure 3.

With this new branch scheme, we can significantly reduce the number of redundant and non-optimal nodes generated in the branch and bound tree.

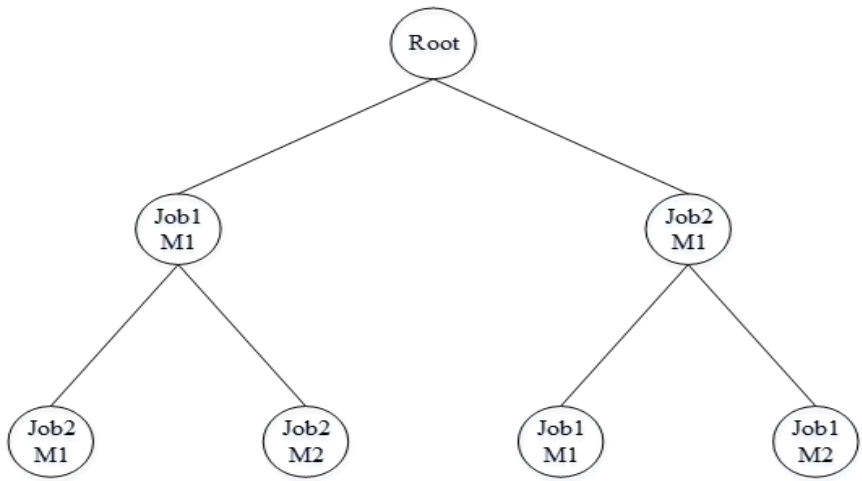

Figure 3

*Lower bound for the optimal objective value*

In this section, we use a heuristic similar to the one introduced by M.Ranjbar[1]. In the heuristic, an initial solution leading to a lower bound is generated and whenever a better solution is obtained during the search process, the lower bound is updated. The heuristic procedure is described in the following pseudocode.

1. Assign $J_i$ to $M_i$ for i=1,…, m. Let k=m and go to step 5;

2. Let $AD\mu_{M_i}$ be the workload of $M_i$ up to current and compute the average load per machine $\bar{\mu} = \dfrac{\sum_i AD\mu_{M_i}}{m}$;

3. For each $M_i$, let $\rho_i = \dfrac{AD\mu_{M_i} - \bar{\mu}}{\sigma_{M_i}}$;

4. Select $J_k$ and assign it to the machine with the smallest value of $\rho$;

5. If k≥n, then stop; otherwise, let k=k+1 and go to step 2.

*Upper bounds for the optimal objective value*

In this section, we introduce three upper bounds for the optimal objective value of the proposed branch and bound algorithm.

**Upper bound 1**

Thorem1: We use $\pi_{u1}$ to represent the objective value of the partial sequence up to current, then $\pi_{u1}$ is an upper bound.

**Proof.** The proof is Straightforward. In this problem, adding a job to a partial sequence will lead to an increase of the makespan and a decrease of the objective value. So the probability of the partial sequence $\pi_{u1}$ is an upper bound of the objective value.

**Upper bound 2**

Thorem 2: The workload of the machines until current is the summation of the total mean time and the setup time of the assigned jobs, We use $m*\delta$ as the X in the standard normal distribution, and use the minimum setup time as the setup time of the unassigned jobs. So the new mean in the standard normal distribution becomes

$$\sum_j MachineLoad_j + \sum \mu_{unassignedJobs} + unassignedjobs * Setup_{min}$$

the standard deviation is still $\sqrt{\sum_i \sigma_i^2}$, so we can get the upper bound of the objective value by the following formula.

$$\phi(\frac{m\delta - \left(\sum_j MachineLoad_j + \sum \mu_{unassignedJobs} + unassignedjob * Setup_{min}\right)}{\sqrt{\sum_i \sigma_i^2}})$$

**Proof.** During the process, the completion time of each machine must be less than or equal to the makespan, we use C as the makespan, $C_i$ as the completion time of the $i$th machine, then we can get $mC \geq \sum_{i=1}^{n} C_i$, then we use $\sum_j MachineLoad_j$ as the workload of all the machines until current, we can know the mean time of the unassigned jobs, and we use $Setup_{min}$ as the setup time of the unassigned jobs. Then we can get

$$mC \geq \sum_j MachineLoad_j + \sum \mu_{unassignedJobs} + (n - assignedjobs)Setup_{min}$$

Now we replace C with $\delta$, and then we can get upper bound of the objective value.

$$upperbound2 = \phi(\frac{m\delta - \left(\sum_j MachineLoad_j + \sum \mu_{unassignedJobs} + unassignedjob * Setup_{min}\right)}{\sqrt{\sum_i \sigma_i^2}})$$

**Upper bound 3**

To find the upper bound 3, we introduce a heuristic process. If this node is in the *i*th level, First we assumed (n-s) identical virtual jobs, the mean processing time for the job is $\mu'$ which is equal to the minimum mean processing time of the unassigned jobs and its variance is $\sigma'^2$ which is equal to the maximum variance of the unassigned jobs, and the setup time is equal to the minimum of all the setup times.

Theorem 3: Then we use the following heuristic rule to get the upper bound of the objective value.

a) Compute the current probability of each machine's completion time not exceed the due date;

b) Assign a virtual job on the machine which has the largest probability;

c) Repeat step a and b, until the virtual jobs are all assigned.

d) Compute the current objective value, and this is the upper bound of the objective value.

**Proof.** We can use the mathematical induction to prove this upper bound.

1) When there is only one job is not assigned, we put a virtual job on the machine which has the largest probability, and then the objective value is the upper bound of the objective value.

2) When there are (n-1) jobs to be assigned, we assign the virtual jobs following the heuristic rule, we assume the current objective value is the upper bound of the objective value.

3) If there are n jobs to be assigned, we have known that when there are (n-1) jobs left, we can get the upper bound of the objective value according to step 2, then on the basis of step 1, we put the last virtual job on the machine which has the largest probability, and the current objective value must be the upper bound of the objective value.

So upper bound 3 is the upper bound of the objective value.

*Dominance rules*

Dominance 1.

Theorem 4: In a node of the branch tree, if the upper bound of the objective value is less than the known overall lower bound, then the node can be cut off.

This theorem is straightforward, and the proof is omitted in this paper.

To find the dominance 2, we refer to Chenjie's work[15] on deterministic parallel machine scheduling problem, and transform it into the stochastic model.

Dominance 2.

Theorem 5: Assume that partial solution P1 and partial solution P2 exist, and they both have assigned the same jobs to the machines. If the workload of all the machines in P1 is no larger than that in P2 and the objective value of P1 up until current is no larger than P2, then P1 dominates P2, and P2 should be eliminated in the B&B process.

**Proof.** Because we have assigned same jobs on the machines now, so the mean of the total processing time and the variance of the total processing time on the machines is the same, if the workload of all the machines in P1 is no larger than that in P2, it must be the setup time in P1 until current is smaller, and because the objective value of P1 until current is better than P2, if we continue assign jobs, the performance of P1 must be better than P2.

# 5. Computational results and the discussion

*Computational setup*

We consider five values for the number of jobs, n=10, 12, 14, 16 and 18, and four values for the number of machines, m=2, 3, 4 and 5. For each job j, the value μ$_j$ is drawn from the normal distribution N[20,9], and the variance is drawn from uniform distribution on interval (0, 0.1 $\mu_i^2$ ), the value of setup time is drawn from uniform distribution on interval [5,10]. The due date is chosen as d= ($\bar{\mu}+\bar{s}+\sqrt{m\bar{\sigma}}$ -5). During the experimentation, we use the best first branch and bound algorithm.

*Performance of the technics*

In this section, we'll first verify the effectiveness of the technic we introduced in the branch and bound algorithm, then we will compare the results between deterministic method and robust method in solving the uncertain problems, and finally we'll present the results which get in reasonable time. If we can't get the results without the new technic in a reasonable time, we use "---" instead.

We use the following GAP formula to indicate the efficiency of the new branch scheme and we represent the CPU time efficiency improvement in Table 1 and the search nodes efficiency improvement in Table 2.

$$GAP = \frac{performance\ without\ new\ technic - performance\ with\ new\ technic}{performance\ without\ new\ technic} * 100\%$$

**Table 1**

Average CPU time gap for the branch scheme.

|   | n |   |   |   |
|---|---|---|---|---|
|   | 10 | 12 | 14 | 16 |
| m |   |   |   |   |
| 2 | 0 | 100 | 588 | --- |
| 3 | 0 | 300 | 1228 | --- |
| 4 | 0 | 0 | 507 | --- |
| 5 | 0 | 142 | 625 | --- |

**Table 2**

Average search nodes gap for the branch scheme.

| | n | | | |
|---|---|---|---|---|
| | 10 | 12 | 14 | 16 |
| m | | | | |
| 2 | 135 | 174 | 291 | --- |
| 3 | 167 | 228 | 286 | --- |
| 4 | 134 | 8 | 270 | --- |
| 5 | 121 | 324 | 612 | --- |

We also use the same formula to test and verify the efficiency of the upper bounds. We get the CPU time improvement in table 3 and the number of search nodes improvement in table 4.

**Table 3**

Average CPU time gap for the upper bounds.

| (m, n) | 2,12 | 2,14 | 3,12 | 3,14 | 4,12 | 5，14 |
|---|---|---|---|---|---|---|
| Only UB1 | 33 | 2 | 50 | 119 | 100 | 64 |
| Only UB2 | 23 | 43 | 40 | 0 | 14 | 0 |
| Only UB3 | 100 | 560 | 100 | 483 | 500 | 174 |
| All Bounds | 300 | 633 | 200 | 900 | 543 | 254 |

**Table 4**

Average search nodes gap for the upper bounds.

| (m, n) | 2,12 | 2,14 | 3,12 | 3,14 | 4,12 | 5，14 |
|---|---|---|---|---|---|---|
| Only UB1 | 3 | 11 | 21 | 8 | 29 | 13 |
| Only UB2 | 20 | 20 | 18 | 14 | 18 | 33 |
| Only UB2 | 398 | 465 | 285 | 453 | 514 | 781 |
| All Bounds | 418 | 465 | 328 | 465 | 531 | 821 |

*Impact of the dominance rule*

In the same way, we test the efficiency of the dominance rule and represent the CPU time improvement in table 5 and the number of search nodes improvement in table 6.

**Table 5**

Average CPU time gap for the dominance rule.

| CPU time | n | | | |
|---|---|---|---|---|
| | 10 | 12 | 14 | 16 |
| m | | | | |
| 2 | 0 | 300 | 44 | 74 |
| 3 | 0 | 900 | 1228 | --- |
| 4 | 100 | 300 | 507 | --- |
| 5 | 100 | 100 | 1200 | --- |

**Table 6**

Average search nodes gap for the dominance rule.

| Nodes | n | | | |
|---|---|---|---|---|
| | 10 | 12 | 14 | 16 |
| m | | | | |
| 2 | 350 | 148 | 304 | 61 |
| 3 | 747 | 437 | 731 | --- |
| 4 | 84 | 58 | 623 | --- |
| 5 | 642 | 100 | 1014 | --- |

*Performance of the initial solution*

We use the following GAP formula to indicate the performance of the initial lower bound, the smaller the GAP, it means the better of the performance of the initial lower bound.

$$GAP = \frac{initial\ lower\ bound - optimal\ objective\ value}{optimal\ objective\ value} * 100\%$$

**Table 7**

Average gap for the initial low bound.

|   | n |   |   |   |
|---|---|---|---|---|
|   | 10 | 12 | 14 | 16 |
| m |   |   |   |   |
| 2 | 137 | 28 | 76 | 17 |
| 3 | 180 | 43 | 0 | 20 |
| 4 | 15 | 14 | 35 | 0 |
| 5 | 21 | 55 | 71 | 10 |

*Comparison between deterministic methods and robust methods*

We solve the stochastic problem with the deterministic method and the robust methods respectively, and use the objective value improvement to verify the efficiency of the robust method.

$$GAP = \frac{obj\ of\ robust\ method - obj\ of\ deterministic\ method}{obj\ of\ deterministic\ method} * 100\%$$

**Table 8**

Average gap between the deterministic methods and robust methods.

|   | n |   |   |   |
|---|---|---|---|---|
|   | 10 | 12 | 14 | 16 |
| m |   |   |   |   |
| 2 | 11 | 29 | 2 | 1 |
| 3 | 68 | 23 | 32 | 26 |
| 4 | 31 | 12 | 782 | 7 |
| 5 | 227 | 72 | 144 | 24 |

*Experiment results*

We present the results of the robust parallel machine scheduling with sequence dependent setup time solved by the branch and bound algorithm within 7200s (2 hours) in Table 9.

**Table 9**

Summary of the results for branch and bound algorithm.

| n  | m | CPU Time | Search Nodes |
|----|---|----------|--------------|
| 10 | 2 | 0        | 1154         |
|    | 3 | 0        | 3344         |
|    | 4 | 0        | 3723         |
|    | 5 | 0        | 1454         |
| 12 | 2 | 1        | 14055        |
|    | 3 | 1        | 14307        |
|    | 4 | 1        | 10829        |
|    | 5 | 5        | 17130        |
| 14 | 2 | 22       | 67659        |
|    | 3 | 28       | 70121        |
|    | 4 | 29       | 55325        |
|    | 5 | 46       | 77215        |
| 16 | 2 | 292      | 247161       |
|    | 3 | 505      | 326384       |
|    | 4 | 616      | 308284       |
|    | 5 | 741      | 309851       |
| 18 | 2 | 1815     | 394781       |
|    | 3 | 2003     | 424084       |
|    | 4 | 2564     | 733581       |
|    | 5 | 7120     | 861719       |

# 6. Conclusions and suggestions

In this paper, we study the robust parallel machine scheduling problem with sequence dependent setup time, which aims to maximize the customer service level in an identical parallel machine environment. A branch and bound algorithm was proposed to solve this problem. Moreover, a new branch scheme, three upper bounds, a lower bound and two dominance rules were developed for the branch and bound algorithm. Computational results revealed that the technic proposed in the branch and bound algorithm is effective, it also revealed the satisfactory performance of the branch and bound algorithm and indicated that it has a good result in solving the robust parallel machine scheduling problem with sequence dependent setup time.

For future work, it is recommended to solve this problem in the more complicated environment such as flow shop environment and job shop environment. Also, the other exact, heuristic and meta-heuristic methods for the robust parallel machine scheduling problem with sequence dependent setup time may also be interesting topics. Finally, it is recommended to try distribution functions other than the normal distribution for this problem.